\documentclass[12pt,oneside,english]{amsart}
\usepackage[T1]{fontenc}
\usepackage[latin9]{inputenc}
\usepackage{geometry}
\usepackage{amssymb}

\makeatletter
\numberwithin{equation}{section} %% Comment out for sequentially-numbered
\numberwithin{figure}{section} %% Comment out for sequentially-numbered
  \@ifundefined{theoremstyle}{\usepackage{amsthm}}{}
  \theoremstyle{plain}
  \newtheorem{thm}{Theorem}[section]
  \theoremstyle{remark}
  \newtheorem{acknowledgement}[thm]{Acknowledgement}
  \theoremstyle{definition}
  \newtheorem{defn}[thm]{Definition}
  \theoremstyle{remark}
  \newtheorem*{note*}{Note}
  \theoremstyle{plain}
  \newtheorem{prop}[thm]{Proposition}
  \theoremstyle{plain}
  \newtheorem{cor}[thm]{Corollary}
  \theoremstyle{plain}
  \newtheorem{lem}[thm]{Lemma}
  \theoremstyle{remark}
  \newtheorem*{notation*}{Notation}

\@ifundefined{definecolor}
 {\usepackage{color}}{}

\usepackage{diagrams}
\usepackage{pinlabel}

\theoremstyle{definition}

\newtheorem*{EmbedThm}{Lagrangian Embedding Theorem}

\newcommand{\R}{\mathbb R}

\newcommand{\N}{\mathbb N}

\newcommand{\micro}{\bf mic}
\newcommand{\ext}{\bf ext}

\newcommand{\vect}{\bf Vect}

\newcommand{\alg}{\bf Alg}

\newcommand{\Sympl}{\bf Sympl}
\newcommand{\ExtSympl}{{\Sympl}^{\ext}}
\newcommand{\MicroSympl}{{\Sympl}_{\micro}}
\newcommand{\ExtMicroSympl}{{\Sympl}_{\micro}^{\ext}}

\newcommand{\ExtMicSym}{\ExtMicroSympl}

\newcommand{\Core}{{\bf core}}

\newcommand{\GR}{{\bf Gr}}

\newcommand{\neutral}{\textbf{E}}
\newcommand{\E}{\textbf{E}}
\newcommand{\neutralmorph}{\textbf{e}}
\newcommand{\e}{\textbf{e}}

\newcommand{\Endop}{\mathcal{END}}

\newcommand{\Mon}[1]{{\alg}(#1)}

\newcommand{\Cot}{{T}^*}
\newcommand{\Tan}{{T}}

\newcommand{\id}{\operatorname{id}}

\newcommand{\graph}{\operatorname{gr}}

\newcommand{\cH}{\mathcal{H}}

\newcommand{\Red}{\operatorname{Red}}
\newcommand{\pr}{\operatorname{pr}}

\makeatother

\makeatother

\usepackage{babel}

\begin{document}

\title{Symplectic Microgeometry I:\\
Micromorphisms}

\author{Alberto S. Cattaneo}
\address{Institut f\"ur Mathematik\\
Universit\"at Z\"urich--Irchel\\
Winterthurerstrasse 190, CH-8057 Z\"urich\\
Switzerland \\
{\it email:} alberto.cattaneo@math.uzh.ch}

\author{Benoit Dherin}
\address{Department of Mathematics\\
Utrecht University\\
Budapestlaan 6, 3584 CD Utrecht\\
The Netherlands\\
{\it email:} b.r.u.dherin@uu.nl }

\author{Alan Weinstein}
\address{Department of Mathematics\\
University of California\\
Berkeley, CA 94720-3840
USA\\ 
{\it email:} alanw@math.berkeley.edu}

\begin{abstract}
We introduce the notion of symplectic microfolds and symplectic micromorphisms
between them. They form a monoidal category, which is a version of
the {}``category'' of symplectic manifolds and canonical relations
obtained by localizing them around lagrangian submanifolds in the
spirit of Milnor's microbundles. 
\end{abstract}
\maketitle
\tableofcontents{}

\section{Introduction}

There is a category $\Sympl$ whose objects are finite-dimensional
symplectic manifolds $(M,\omega)$ and whose morphisms are symplectomorphisms
$\Psi:(M,\omega_{M})\rightarrow(N,\omega_{N})$. In attempting to
understand the quantization procedure of physicists from a mathematical
perspective, one may think of it as a functor from this symplectic
category, where classical mechanics takes place, into the category
of Hilbert spaces and unitary operators, which is the realm of quantum
mechanics.

It is well known that the category $\Sympl$ is too large, since there
are {}``no-go'' theorems that show that the group of all symplectomorphisms
on $(M,\omega)$ does not act in a physically meaningful way on a
corresponding Hilbert space. One standard remedy for this is to replace
$\Sympl$ with a smaller category, replacing the symplectomorphism
groups with certain finite-dimensional subgroups. Another is to replace
the Hilbert spaces and operators by objects depending on a formal
parameter.

But there is also a sense in which the category $\Sympl$ is too {\em
small}, since it does not contain morphisms corresponding to operators
such as projectors and the self-adjoint (or skew-adjoint) operators
that play the role of observables in quantum mechanics, nor can it
encode the algebra structure itself on the space of observables. (This
collection of observables is not actually a Hilbert space, but certain
sets of operators do carry a vector space structure, with the inner
product associated to the Hilbert-Schmidt norm.)

To enlarge the symplectic category, we look at the {}``dictionary''
of quantization, following, for example, \cite{weinstein1996}. In
this dictionary, the cartesian product of symplectic manifolds corresponds
to the tensor product of Hilbert spaces, and replacing a symplectic
manifold $(M,\omega)$ with $(M,-\omega)$ (which we denote by $\overline{M}$
when we omit the symplectic structure from the notation for a given
symplectic manifold) corresponds to replacing a Hilbert space $\cH$
by its conjugate, or dual, space $\cH^{*}$. Thus, if symplectic manifolds
$M_{1}$ and $M_{2}$ correspond to Hilbert spaces $\cH_{1}$ and
$\cH_{2}$, the product $\overline{M}\times N$ corresponds to $\cH_{1}^{*}\otimes\cH_{2}$,
which, with a suitable definition of the tensor product, is a space
$L(\cH_{1},\cH_{2})$ of linear operators from $\cH_{1}$ to $\cH_{2}$.

Another entry in the dictionary says that lagrangian submanifolds
(perhaps carrying half-densities) in symplectic manifolds correspond
to vectors or lines in Hilbert space. Combining this idea with the
one in the previous paragraph, we conclude that lagrangian submanifolds
in $\overline{M}\times N$ should correspond to linear operators from
$\cH_{1}$ to $\cH_{2}$.This suggests that, if the space of observables
$\cH$ for a quantum system corresponds to a symplectic manifold $M$,
then the algebra structure on $\cH$ should be given by a lagrangian
submanifold $\mu$ in $\overline{M}\times\overline{M}\times M.$ The
algebra axioms of unitality and associativity should be encoded by
monoidal properties of $\mu$ in an extended symplectic category,
$\ExtSympl$, where the morphisms from $M$ to $N$ are the canonical
relations; i.e., all the lagrangian submanifolds of $\overline{M}\times N$
(not just those which are the graphs of symplectomorphisms) and where
the morphism composition is the usual composition of relations. However,
a problem immediately occurs: the composition of canonical relations
may yield relations that are not submanifolds any more, and thus are
not canonical relations! $\ExtSympl$ is then not a true category,
as the morphisms cannot always be composed. It is rather awkward to
speak about a quantization functor in this context.

There have already been several approaches to remedy this defect.
One approach, developed by Guillemin and Sternberg in \cite{GS1979}
(see \cite{GH2007} for a recent version), is to consider only symplectic
vector spaces and linear canonical relations. Another, suggested by
Wehrheim and Woodward in \cite{WW2007}, is to enlarge the category
still further by allowing arbitrary {}``formal'' products of canonical
relations and equating them to actual products when the latter exist
as manifolds.

In this paper, we take yet another approach. We construct a version
of the extended symplectic {}``category,'' which is a true category,
by localizing it around lagrangian submanifolds. Its objects, called
symplectic microfolds in the spirit of Milnor's microbundles (\cite{milnor1964}),
are equivalence classes $[M,A]$ of pairs consisting of a symplectic
manifold $M$ and a lagrangian submanifold $A\subset M$, called the
core. The equivalence reflects the fact that these objects really
describe the geometry of a neighborhood $-$ or a {}``micro'' neighborhood
$-$ of $A$ in $M$. 

In this {}``micro'' setting, there is also a notion of canonical
{}``micro'' relations between two symplectic microfolds: They are
lagrangian submicrofolds $[L,C]$ of the symplectic microfold product
$[\overline{M},A]\times[N,B]$. Their composition is generally as
ill behaved as it is for regular canonical relations. 

One of the main points of this paper is to identify a certain subset
of canonical {}``micro'' relations satisfying a new transversality
condition which ensures the composition is always well defined. We
consider these transverse {}``micro'' canonical relations as morphisms
between symplectic microfolds; in this way, we obtain a new symmetric
monoidal category: the extended microsymplectic category. 

The extended symplectic {}``category'' has been used as a sort of
heuristic guideline in an attempt to quantize Poisson manifolds (see
\cite{WW1991}) in a geometric way. These attempts have been only
partially successful due in part to the existence of nonintegrable
Poisson manifolds (hence restricting the class of Poisson manifold
one can quantize), as well as to the ill-defined composition of canonical
relations (limiting thus the functorial properties of these geometric
quantization methods). The replacement of the extended symplectic
{}``category'' by its {}``micro'' version provides new ways of
dealing with both issues. 

This paper lays the foundation for a series of work that revolves
around two main themes: the categorification of Poisson geometry and
its functorial quantization as explained below.

\subsection*{Categorification}

Since $\ExtMicroSympl$ is a monoidal category, it is natural to consider
its category of algebras $\Mon{\ExtMicroSympl}$. The main statement
we are aiming at here is the equivalence between this latter category
and the category of Poisson manifolds and Poisson maps. Future research
directions will include the study of a weakened version of $\Mon{\ExtMicroSympl}$
whose algebra maps are replaced by bimodules. This should correspond
to a {}``micro'' Morita theory for Poisson manifolds.

\subsection*{Quantization}

Our second line of work will focus on constructing a monoidal functor
from the extended microsymplectic category (enhanced with half-density
germs on the morphisms) to the category of vector spaces. Since monoidal
functors between monoidal categories induce functors between their
respective categories of monoid objects, we obtain in this way a {}``quantization''
functor from the category of Poisson manifolds to the category of
algebras in $\vect$. 

\begin{acknowledgement}
We thank Giovanni Felder, Domenico Fiorenza, Jim Stasheff and Chenchang
Zhu. A.S.C. acknowledges partial support of SNF Grant 200020-121640/1,
of the European Union through the FP6 Marie Curie RTN ENIGMA (contract
number MRTN-CT-2004-5652), and of the European Science Foundation
through the MISGAM program. B. D. acknowledges partial support from
SNF Grant PA002-113136 and from the Netherlands Organisation for Scientific
Research (NWO), and thanks Wendy L. Taylor for proofreading the manuscript.
A.W. acknowledges partial support from NSF grant DMS-0707137. 
\end{acknowledgement}

\section{Symplectic microfolds\label{def:MfdPairs}}

A $\mathbf{local}$ $\mathbf{manifold}$ $\mathbf{pair}$ $(M,A)$
consists of a manifold $M$ and a submanifold $A\subset M$, called
the $\mathbf{core}$. Two manifold pairs $(M,A)$ and $(N,B)$ are
said to be equivalent if $A=B$ and if there is a third manifold pair
$(U,A)$ such that $U$ is an open subset in both $M$ and $N$ simultaneously.
A map between local pairs is a smooth map from $M$ to $N$ that sends
$A$ to $B$. Note that we require equality of neighborhoods and not
merely diffeomorphism.

\begin{defn}
A $\mathbf{microfold}$ is an equivalence class of a local pairs $(M,A)$.
We denote these equivalence classes either by $[M,A]$ or by $([M],A)$.
Sometimes, $[M]$ will be referred to as a $\mathbf{manifold}$ $\mathbf{germ}$
around $A$. 
\end{defn}
We define an equivalence relation on the maps of local pairs that
send a representative of $[M,A]$ to a representative of $[N,B]$
by declaring two such maps equivalent if there is a common neighborhood
of $A$ where they coincide. The equivalence classes, written as\begin{eqnarray*}
[\Psi]:[M,A] & \longrightarrow & [N,B],\end{eqnarray*}
 are the maps between microfolds. We say that $[\Psi]$ is a germ
above $\phi:A\rightarrow B$ if, for $\Psi\in[\Psi],$ we have that
$\Psi_{|A}=\phi$.

A $\mathbf{submicrofold}$ of a microfold $[M,A]$ is a microfold
$[N,B]$ such that $N\subset M$ and $B\subset A$. We define the
graph of a microfold map $[\Psi]$ as the submicrofold\begin{eqnarray*}
\graph[\Psi] & := & \big([\graph\Psi],\graph\Psi_{|A}\big)\end{eqnarray*}
of the product microfold \begin{eqnarray*}
[M,A]\times[N,B] & := & [M\times N,A\times B].\end{eqnarray*}

Microfolds and microfold maps form a category. Fibered products of
microfolds are defined in the obvious way.

There are micro counterparts of symplectic manifolds.

\begin{defn}
A $\mathbf{symplectic}$ $\mathbf{microfold}$ is a microfold $[M,L]$
where $M$ is a symplectic manifold and $L\subset M$ a lagrangian
submanifold. We call $\mathbf{cotangent}$ $\mathbf{microbundles}$
the symplectic microfolds of the form $[\Cot M,M].$ 
\end{defn}
\begin{note*}
We will write $Z_{E}$ to denote the zero section of a vector bundle
$E\rightarrow M$. In the previous definition, we abused notation
by writing $[\Cot M,M]$ instead of $[\Cot M,Z_{\Cot M}]$. 
\end{note*}
A symplectomorphism between symplectic microfolds is a microfold map
for which there is a representative that is a symplectomorphism.

\begin{note*}
To avoid an explosion in the use of the prefix {}``micro'', we will
keep the usual {}``manifold terminology'' when available and assume
that we are talking about the {}``micro'' version when microfolds
are around and no confusion is possible. For examples, we choose to
use {}``symplectomorphism'' instead of {}``microsymplectomorphism'',
and so on. 
\end{note*}
Symplectic microfolds and their symplectomorphisms form a category,
which we denote by $\MicroSympl$.

Many special submanifolds of symplectic geometry have their corresponding
{}``micro'' versions.

\begin{defn}
A microsubmanifold $[S,X]$ of a symplectic microfold $[M,L]$ will
be called $\mathbf{isotropic}$, $\mathbf{lagrangian}$, or $\mathbf{coisotropic}$
if there are representatives of $[S]$ and $[M]$ which are isotropic,
lagrangian or coisotropic, respectively. 
\end{defn}
The language of microfolds is useful to express local geometric properties;
that is, properties that are true for all neighborhoods of some submanifold.
For instance, the lagrangian embedding theorem can be phrased as follows.

\begin{EmbedThm}For any symplectic microfold $[M,L]$, there exists
a symplectomorphism \begin{eqnarray}
{}[\Psi_{M,L}]:[M,L] & \longrightarrow & [\Cot L,L],\label{Thm: lagrangina embedding}\end{eqnarray}
 above the identity. \end{EmbedThm}Actually, this theorem was first
stated and proved using the language of local manifold pairs and their
equivalences (see \cite{Weinstein1971}).

An important notion in symplectic geometry is that of a $\textbf{canonical}$
$\mathbf{relation}$; i.e., a lagrangian submanifold $V\subset\overline{M}\times N$,
where $(M,\omega_{M})$ and $(N,\omega_{N})$ are symplectic manifolds
and $\overline{M}$ is the symplectic manifold $(M,-\omega_{M})$.
Canonical relations are usually thought of as {}``generalized symplectomorphisms''
and written as $V:M\rightarrow N$ instead of $V\subset\overline{M}\times N$.
The rationale behind this is twofold: the graph of a symplectomorphism
is a canonical relation, and it is formally possible to extend the
composition of symplectomorphisms to canonical relations. Namely,
the composition of $V\subset\overline{M}\times N$ and $W\subset\overline{N}\times P$
is the subset of $\overline{M}\times P$ defined by\begin{eqnarray*}
W\circ V & = & \Red\Big(\big(V\times W\big)\cap\big(M\times\Delta_{N}\times P\big)\Big)\\
 & = & \Red\big(V\times_{N}W\big),\end{eqnarray*}
 where $\Red$ is the {}``reduction'' map that projects $(m,n,n,p)$
to $(m,p)$. The major issue here is that composition of canonical
relations is generally ill defined: $W\circ V$ may fail to be a submanifold,
although, when it is, it is a lagrangian one%
\footnote{This is a special instance of symplectic reduction. Namely, the submanifolds
$M\times\Delta_{N}\times P$ and $V\times W$ are respectively coisotropic
and lagrangian in $\overline{M}\times N\times\overline{N}\times P$.
The composition $W\circ V$ is exactly the quotient of $V\times W$
by the characteristic foliation of $M\times\Delta_{N}\times P$. This
ensures that $W\circ V$ is a lagrangian submanifold whenever it is
a submanifold. %
}. There is a well known criterion which limits the wildness of the
composition and which we will need later.

\begin{thm}
\label{thm:clean intersection}The composition $W\circ V$ of the
canonical relations $V$ and $W$ as above is an immersed lagrangian
submanifold of $\overline{M}\times P$ if the submanifolds $V\times W$
and $M\times\Delta_{N}\times P$ intersect cleanly. 
\end{thm}
Nevertheless, it is standard to think of symplectic manifolds and
canonical relations as a category. It is called the $\mathbf{extended}$
$\mathbf{symplectic}$ {}``$\mathbf{category}$''%
\footnote{The quotes are there as a reminder that it is not really a category.%
} and will be denoted by $\ExtSympl$. Many constructions in $\Sympl$
extend to $\ExtSympl$. For instance, we define the image of a point
$x\in M$ by a canonical relation $V\subset\overline{M}\times N$
as the subset of $N$ given by\begin{eqnarray}
V(x) & := & \pi_{N}\Big(V\cap\big(\{x\}\times N\big)\Big),\label{def: Image trough relation}\end{eqnarray}
 where $\pi_{N}$ is the projection on the second factor of $M\times N$.
The tangent relation $\Tan V:\Tan M\rightarrow\Tan N$ to a canonical
relation $V:M\rightarrow N$ as the subset $\Tan V$ of $\Tan M\times\Tan N$
given by the set of tangent vectors to $V$. 

The notion of canonical relation between symplectic manifolds can
be transported to symplectic microfolds.

\begin{defn}
A $\textbf{canonical relation}$ $([V],K)$ between the symplectic
microfolds $[M,A]$ and $[N,B]$ is a lagrangian submicrofold $([V],K)$
of $[\overline{M}\times N,A\times B]$. 
\end{defn}
\begin{note*}
We will often prefer the notation $([V],K)$ for canonical relations
and reserve the notation $[M,A]$ for symplectic microfolds in order
to distinguish between objects and morphisms. We will also use the
arrow notation\begin{eqnarray*}
\big([V],K\big):[M,A] & \longrightarrow & [N,B]\end{eqnarray*}
 to represent canonical relations between symplectic microfolds. We
interpret the core $K$ as a submanifold of $B\times A$ and not as
a submanifold of $A\times B$. In other words, we decide to regard
the core as a generalized morphism from $B$ to $A$. The reason for
this contraintuitive interpretation will become apparent later on. 
\end{note*}
The composition of canonical relations in the microworld \[
\begin{array}{ccccc}
[M,A] & \overset{([V],K)}{\longrightarrow} & [N,B] & \overset{([W],L)}{\longrightarrow} & [P,C]\end{array}\]
 is given by the binary relation composition of their {}``components''
\begin{eqnarray}
([W],L)\circ([V],K) & := & ([W\circ V],K\circ L).\label{eq:composition}\end{eqnarray}
 At this point, we still have the same kind of ill defined composition
for canonical relations between symplectic microfolds as we had for
symplectic manifolds. However, we may consider a special type of canonical
relations between symplectic microfolds that always compose well.
This is what we do next.

\section{Symplectic micromorphisms\label{sub:Generalized-symplectomorphisms}}

\subsection{Definitions}

Our starting point is the cotangent lift $\Cot\phi$ of a diffeomorphism
$\phi$:\begin{diagram}
   \Cot A & \rTo{\Cot \phi} & \Cot B \\
   \dTo{\pi_A}   &               & \dTo_{\pi_B}\\
   A       & \lTo_{\phi}   &     B
\end{diagram}It induces a canonical relation of the form\begin{eqnarray*}
\big([\graph\Cot\phi],\graph\phi\big):[\Cot A,A] & \longrightarrow & [\Cot B,B].\end{eqnarray*}
The fact that the canonical relation $\big([\graph\Cot\phi],\graph\phi\big)$
comes from a map implies the identities%
\footnote{Note that, in general, the graph of a map $f:X\rightarrow Y$, seen
as a generalized morphism from $Y$ to $X$, satisfies $\big(\graph f\big)(y)=f^{-1}(y),$
for $y\in Y$. %
} \begin{eqnarray}
\Big(\graph\Cot\phi\Big)(x) & = & \phi^{-1}(x),\label{eq:Cot 1}\\
\Big(\Tan\graph\Cot\phi\Big)(v) & = & \big(\Tan\phi\big)^{-1}(v),\label{eq:Cot 2}\end{eqnarray}
for all $x$ in the zero section of $\Cot A$ and all tangent vectors
$v$ to the zero section of $\Cot A$. Obviously, canonical relations
coming from cotangent lifts compose well. It turns out that the identities
\eqref{eq:Cot 1} and \eqref{eq:Cot 2} are the key to this nice composability.
We therefore make the following definition.

\begin{defn}
A $\mathbf{symplectic}$ $\mathbf{micromorphism}$ is a canonical
relation of the form\begin{eqnarray*}
\big([V],\graph\phi\big):[M,A] & \longrightarrow & [N,B],\end{eqnarray*}
 where $\phi$ is a smooth map from $B$ to $A$ and such that there
exists a representative $V\in[V]$ for which\begin{eqnarray}
V(a) & = & \phi^{-1}(a),\quad\textrm{for all }a\in A,\label{eq:CondI bis}\\
TV(v) & = & \big(T\phi\big)^{-1}(v),\quad\textrm{for all }v\in T_{a}A.\label{eq:Cond II bis}\end{eqnarray}
where $V(a)$ is the image of $a$ under the relation $V$ as defined
by \eqref{def: Image trough relation} and $TV(v)$ is the image of
$v$ under the tangent relation $TV$. We will usually write $\big([V],\phi\big)$
instead of $\big([V],\graph\phi\big)$. 
\end{defn}
The next Proposition gives various characterizations of symplectic
micromorphisms. Recall that a submanifold $X$ of a manifold $M$
is transverse to a subbundle $\lambda\rightarrow Y$ of $\Tan M\rightarrow M$
along a submanifold $Z\subset X\cap Y$ if \[
\Tan_{z}X+\lambda_{z}=\Tan_{z}M,\quad z\in Z.\]
In this case we write $X\pitchfork_{Z}\lambda$. In particular, a
submanifold $X$ is transverse to a submanifold $Y$ along $Z\subset X\cap Y$
if $X$ is transverse to $\Tan Y$ along $Z$; we write this $X\pitchfork_{Z}Y$.
A $\textbf{splitting}$ of a symplectic microfold $[M,A]$ is a lagrangian
subbundle $K$ of the tangent bundle $TM\big|_{A}$ of $M$ restricted
to $A$ such that, for all $x\in A$, we have%
\footnote{Note that we will sometimes abuse notation slightly by writing $T_{x}L$
instead of $T_{x}L\oplus0$ and $L_{x}$ instead of $0\oplus L_{x}$.%
} $T_{x}M=T_{x}A\oplus K_{x}$ (i.e. $A$ is transverse to $K$ along
$A$). 

\begin{defn}
A canonical relation of the form\begin{eqnarray*}
\big([V],\phi\big):[M,A] & \longrightarrow & [N,B]\end{eqnarray*}
is said to be $\mathbf{transverse}$ $\mathbf{to}$ $\mathbf{a}$
$\mathbf{splitting}$ $K$ of $[N,B]$ if there is a $V\in[V]$ such
that\[
V\cap(A\times N)=\graph\phi\quad\textrm{and}\quad V\pitchfork_{\graph\phi}\Big(\big(\Tan A\oplus0\big)\times\big(0\oplus K\big)\Big).\]
In this case, we will abuse notation slightly and write\[
[V]\pitchfork_{\graph\phi}\big(\Tan A\times K\big).\]
$ $If $([V],\phi)$ is transverse to all splittings, we will say
that it is $\mathbf{strongly}$ $\mathbf{transverse}$. 
\end{defn}
The following proposition offers alternative descriptions of symplectic
micromorphisms in terms of transverse intersections as pictured below:
\begin{figure}[h!] 
\labellist \small\hair 2pt 
\pinlabel $TA\times K$ at 10 180 
\pinlabel $([V],\phi)$ at 10 150 
\pinlabel $([W],\phi)$ at 10 60 
\pinlabel $A\times N$ at 140 30 
\pinlabel $\graph \phi$ at 140 60 
\pinlabel $B$ at 187 105 
\pinlabel $N$ at 187 160 
\pinlabel $A$ at 247 40 
\pinlabel $M$ at 289 40 
\pinlabel $A\times N$ at 249 160 
\pinlabel $\graph \phi$ at  237 100 
\pinlabel $([V],\phi)$ at 280 65
\endlabellist 
\centering 
\includegraphics[scale=1]{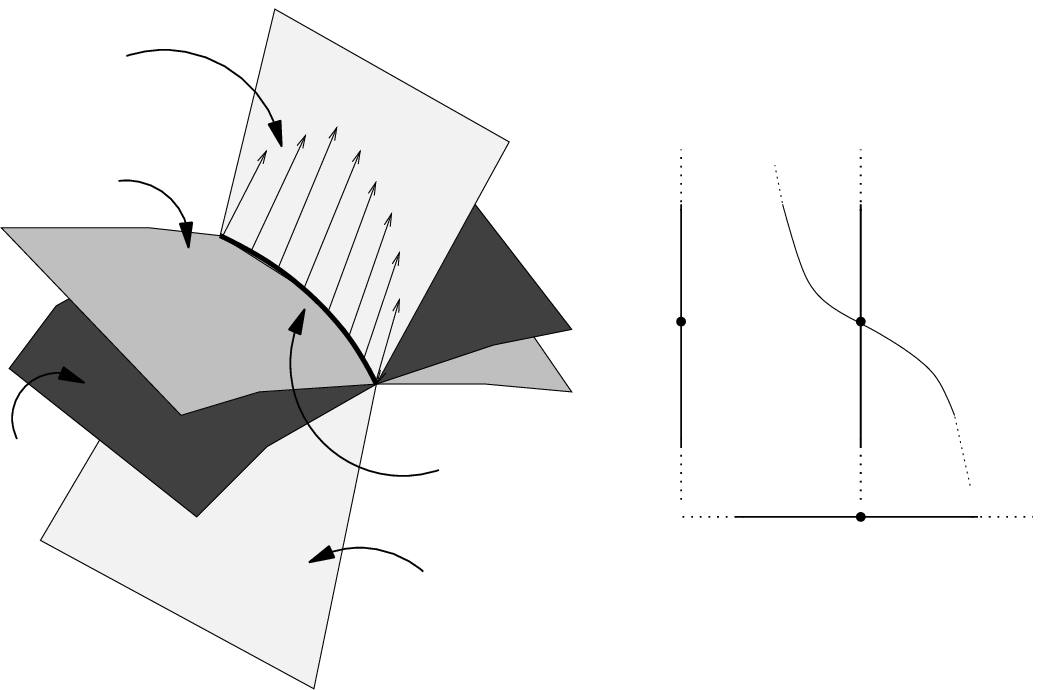}  \end{figure}

\begin{prop}
\label{pro:tranversality}Consider a canonical relation of the form\begin{eqnarray*}
\big([V],\phi\big):[M,A] & \longrightarrow & [N,B].\end{eqnarray*}
Then, the following statements are equivalent:

1. $([V],\phi)$ is a symplectic micromorphism,

2. $([V],\phi)$ is strongly transverse,

3. there is a $V\in[V]$ such that $V\cap(A\times N)=\graph\phi$
is transverse.
\end{prop}
\begin{proof}
$\underline{\textrm{First step}}:\,1\Rightarrow2$. In general, we
have that\begin{eqnarray*}
V\cap(A\times N) & = & \bigcup_{a\in A}\{a\}\times V(a),\end{eqnarray*}
which yields $V\cap(A\times N)=\graph\phi,$ because of \eqref{eq:CondI bis}.
Similarly, one obtains \begin{eqnarray*}
\Tan V\cap(\Tan A\times\Tan N) & = & \graph T\phi\end{eqnarray*}
from \eqref{eq:Cond II bis}. Now, for any splitting $K$ of $[N,B]$,
we have that\[
TA\times K\subset TA\times TN,\]
and, therefore, that\begin{eqnarray*}
\Tan V\cap(\Tan A\times K) & = & \Big(\Tan V\cap(\Tan A\times TN)\Big)\cap(\Tan A\times K)\\
 & = & \Tan\graph\phi\cap(\Tan A\times K),\end{eqnarray*}
along $\graph\phi$. Since $\Tan\graph\phi\subset\Tan A\times\Tan B$
and since $\Tan B\cap K=0$, we have that\[
(v,w)\in\Tan\graph\phi\cap(\Tan A\times K)\]
if and only if $w=0$. By definition, we have that \begin{eqnarray*}
\Tan_{(\phi(b),b)}\graph\phi & = & \big\{(\Tan_{b}\phi(w),w):\, w\in\Tan_{b}B\big\},\end{eqnarray*}
and, therefore, we can conclude that\begin{eqnarray*}
\Tan V\cap(\Tan A\times K) & = & 0\end{eqnarray*}
along $\graph\phi$. A dimension count yields\[
\dim\Tan_{(\phi(b),b)}V+\dim(\Tan_{\phi(b)}A\times K_{b})=\dim\Tan_{(\phi(b),b)}(M\times N),\]
which completes the proof that $\Tan V$ intersects $\Tan A\times K$
transversally along $\graph\phi$. 

$\underline{\textrm{Second step}}:\,2\Rightarrow3$. First, note that
the condition $V\cap(A\times N)=\graph\phi$ is part of both statements.
Let $K$ be a splitting of $[N,B]$. By hypothesis, we have that\begin{eqnarray*}
\Tan V+(\Tan A\times K) & = & \Tan(M\times N)\end{eqnarray*}
along $\graph\phi$, which implies in particular that\begin{eqnarray*}
\Tan V+(\Tan A\times\Tan N) & = & \Tan(M\times N),\end{eqnarray*}
meaning that $A\cap(A\times N)=\graph\phi$ is a transverse intersection. 

$\underline{\textrm{Third step}}:\,3\Rightarrow1$. The fact that
there exists $V\in[V]$ such that\begin{eqnarray*}
V\cap(A\times N) & = & \graph\phi\end{eqnarray*}
 implies \eqref{eq:CondI bis}. Namely, this gives immediately that
\begin{eqnarray*}
V\cap\big(\{a\}\times N\big) & = & \{a\}\times\phi^{-1}(a),\end{eqnarray*}
and, by definition of $V(a)$ (see \eqref{def: Image trough relation}),
we obtain that $V(a)=\phi^{-1}(a)$. Since the intersection of $V$
with $A\times N$ is transverse, it is also clean, i.e., we have that\begin{eqnarray*}
TV\cap(\Tan A\times\Tan N) & = & \graph\Tan\phi.\end{eqnarray*}
Using this equation and the same argument as above, one obtains \eqref{eq:Cond II bis}.
\end{proof}
In Step 1 of the proof of Proposition \ref{pro:tranversality}, we
showed that if a canonical relation $([V],\phi)$ is a symplectic
micromorphism then there exists $V\in[V]$, such that $V$ intersects
$A\times N$ cleanly in $\graph\phi$: \begin{eqnarray}
V\cap(A\times N) & = & \graph\phi\label{eq:CondI}\\
TV\cap(TA\times TN) & = & \graph T\phi.\label{eq:CondII}\end{eqnarray}
In turn, we proved that that conditions \eqref{eq:CondI} and \eqref{eq:CondII}
imply that $([V],\phi)$ is strongly transverse. This makes conditions
\eqref{eq:CondI} and \eqref{eq:CondII} equivalent to $([V],\phi)$
being a symplectic micromorphism. Therefore, if we compare this with
Statement 3 of Proposition \ref{pro:tranversality}, we see that clean
intersection is enough. In Step 2, we see that it is enough to have
a single splitting $K$ of $[N,B]$ transverse to $([V],\phi)$ to
show that it is a symplectic micromorphism. We obtain thus the following,
apparently weaker, version of Proposition \ref{pro:tranversality}. 

\begin{cor}
\label{cor:tranversality}Consider a canonical relation of the form\begin{eqnarray*}
\big([V],\phi\big):[M,A] & \longrightarrow & [N,B].\end{eqnarray*}
Then, the following statements are equivalent:

1. $([V],\phi)$ is a symplectic micromorphism,

2. $([V],\phi)$ is transverse to a splitting of $[N,B]$. 

3. there is a $V\in[V]$ such that $V\cap(A\times N)=\graph\phi$
is clean. 
\end{cor}

\subsection{Examples\label{sub:Examples}}

\subsubsection*{The unit symplectic microfold $\neutral$}

\label{exa: Unit object and morphism} Let us denote by $\E$ the
cotangent bundle of the one point manifold%
\footnote{We may define $\{\star\}$ as the manifold containing only the singleton
$\star=\{\emptyset\}$. %
} $\{\star\}$, which we regard as a symplectic microfold; i.e.,\begin{eqnarray*}
\neutral & := & \Big[\Cot\{\star\},\big\{(0,\star)\big\}\Big].\end{eqnarray*}
 For any microsymplectic manifold $[M,L]$, there is a unique symplectic
micromorphism\begin{eqnarray*}
\e_{[M,L]}:\E & \longrightarrow & [M,L]\end{eqnarray*}
 given by\begin{eqnarray*}
\e_{[M,L]}: & = & \Big(\big[\{(0,\star)\}\times L\big],\pr_{L}\Big),\end{eqnarray*}
 where $\pr_{L}$ is the unique map from $L$ to $\{\star\}$. On
the other hand, symplectic micromorphisms \begin{eqnarray*}
\nu:[M,L] & \longrightarrow & \neutral\end{eqnarray*}
 are in bijection with lagrangian submanifold germs $[V_{x}]$ around
a given point $x\in L$ that are transverse to $L$ at $x$. Namely,
the core map\begin{eqnarray*}
\Core(\nu):\{\star\} & \longrightarrow & L\end{eqnarray*}
 is specified by the image $x\in M$ of the unique point $\star$;
hence\begin{eqnarray*}
\nu & = & \big([V_{x}],\{x\}\big).\end{eqnarray*}
 Conditions \eqref{eq:CondI bis} and \eqref{eq:Cond II bis} read
$V\cap L=\{x\}$ and $T_{x}V\cap T_{x}L=0$; i.e., such that $V_{x}$
and $L$ are transverse.

\subsubsection*{Cotangent lifts}

\label{ex:cotlift} Recall the standard identification\begin{eqnarray*}
\overline{\Cot M}\times\Cot N & \simeq & \Cot(M\times N)\end{eqnarray*}
 via the {}``Schwartz transform'' (see \cite{weinstein1996})\[
\mathcal{S}\big((p_{1},x_{1}),(p_{2},x_{2})\big)=\big((-p_{1},p_{2}),(x_{1},x_{2})\big).\]
For any smooth map $\phi:N\rightarrow M$, the transform of the conormal
bundle \[
N^{*}(\graph\phi)\subset\Cot(M\times N)\]
is a symplectic micromorphism from $[\Cot M,M]$ to $[\Cot N,N]$
given by \begin{eqnarray*}
\Cot\phi & := & \Big(\big[\mathcal{S}^{-1}\big(N^{*}(\graph\phi)\big)\big],\phi\Big).\end{eqnarray*}
We call it the $\mathbf{cotangent}$ $\mathbf{lift}$ of $\phi$ and
denote it again by $\Cot\phi$. Note that, whenever $\phi$ is a diffeomorphism,
we slightly abuse notation since in this case\begin{eqnarray*}
\Cot\phi & = & \Big(\big[\graph\Cot\phi],\phi\Big).\end{eqnarray*}

\subsubsection*{Symplectomorphism germs}

As in the macroworld, the graph\begin{eqnarray*}
\graph[\Psi] & = & \big([\graph\Psi],\Psi_{|B}^{-1}\big)\end{eqnarray*}
of symplectomorphism germ $[\Psi]:[M,A]\rightarrow[N,B]$ is a symplectic
micromorphism: obviously, we have that\begin{eqnarray*}
(\graph\Psi)\cap(A\times N) & = & \graph\Psi_{B}^{-1},\\
(\graph\Tan\Psi)\cap(\Tan A\times\Tan N) & = & \graph\Tan\Psi_{B}^{-1}.\end{eqnarray*}
The following proposition tells us that every symplectic micromorphisms
whose core map is a diffeomorphism is the graph of a symplectomorphism
germ. 

\begin{prop}
\label{prop:germinv} If the core $\phi$ of a symplectic micromorphism\begin{eqnarray*}
([V],\phi):[M,A] & \rightarrow & [N,B]\end{eqnarray*}
 is a diffeomorphism, then there exists a symplectomorphism germ\begin{eqnarray*}
[\Psi]:[M,A] & \rightarrow & [N,B]\end{eqnarray*}
 such that $\graph[\Psi]=([V],\phi)$. 
\end{prop}
\begin{proof}
For each $a\in A$, we have that the intersection\begin{eqnarray*}
V\cap\big(\{a\}\times N\big) & = & \big(a,\phi^{-1}(a)\big)\end{eqnarray*}
is transverse. We see this by counting the dimensions and by remarking
that tangent space intersection \[
\Tan_{(a,\phi^{-1}(a))}V\cap\big(\{0\}\times\Tan_{\phi^{-1}(a)}N\big)\]
is contained in $\graph\Tan\phi^{-1}$, which implies that this intersection
must be zero since\begin{eqnarray*}
\big(\{0\}\times\Tan_{\phi^{-1}(a)}N\big)\cap\graph\Tan\phi^{-1} & = & \{0\}.\end{eqnarray*}
This transverse intersection guarantees that, for each $a\in A$,
there is a neighborhood $K_{a}$ of $\big(a,\phi^{-1}(a)\big)$ in
$M\times N$ and a neighborhood $U_{a}$ of $a$ in $M$ such that
the restriction of the first factor projection,\begin{eqnarray*}
M\times N\supset V\cap K_{a} & \longrightarrow & U_{a}\subset M,\end{eqnarray*}
is a diffeomorphism. Let us denote by $\rho_{a}$ its inverse and
set $\Psi_{a}:=\pi_{N}\circ\rho_{a}$. By construction, we have that\begin{eqnarray*}
\graph\Psi_{a} & = & V\cap K_{a}.\end{eqnarray*}
Since, for two overlapping neighborhoods $U_{a}$ and $U_{a'}$ as
above, the corresponding maps $\Psi_{a}$ and $\Psi_{a'}$ coincide
on $U_{a}\cap U_{a'}$, we obtain a germ\begin{eqnarray*}
[\Psi]:[M,A] & \longrightarrow & [N,B]\end{eqnarray*}
such that $\graph[\Psi]=([V],\phi)$. Because $[V]$ is a lagrangian
submanifold germ, $[\Psi]$ is a symplectomorphism germ. 
\end{proof}

\subsection{Composition}

We now prove that the composition of symplectic micromorphisms\begin{diagram}
   {} [M,A] & \rTo{\big([V],\phi\big)}&  [N,B] & \rTo{\big([W],\psi\big)} &[P,C].
\end{diagram}is always well defined and that\begin{eqnarray*}
\big([W\circ V],\phi\circ\psi\big):[M,A] & \longrightarrow & [P,C]\end{eqnarray*}
 is again a symplectic micromorphism. We first prove a cleanness result
in order to apply Theorem \ref{thm:clean intersection} to our micro
setting.

\begin{lem}
\label{lem:Transvers}Let $\big([V],\phi\big)$ and $\big([W],\psi\big)$
be symplectic micromorphisms as above. For all $V\in[V]$ and $W\in[W]$,
$V\times W$ intersects $M\times\Delta_{N}\times N$ transversally
along $\graph\phi\times_{B}\graph\psi$. 
\end{lem}
\begin{proof}
We need to show that\begin{eqnarray*}
T(V\times W)+T(M\times\Delta_{M}\times P) & = & T(M\times N\times N\times P)\end{eqnarray*}
at all points\begin{eqnarray*}
K(p) & := & \Big(\phi\circ\psi(p),\psi(p),\psi(p),p\Big)\end{eqnarray*}
in $\graph\phi\times_{B}\graph\psi.$ In the symplectic vector space
\[
T_{\phi\circ\psi(p)}\overline{M}\times T_{\psi(p)}N\times T_{\psi(p)}\overline{N}\times T_{p}P\]
 we have that \begin{eqnarray*}
T_{K(p)}(V\times W)^{\perp} & = & T_{K(p)}(V\times W)\\
(T_{\phi\circ\psi(p)}M\times T_{(\psi(p),\psi(p))}\Delta_{N}\times T_{p}P)^{\perp} & = & \{0\}\times T_{(\psi(p),\psi(p))}\Delta_{N}\times\{0\}.\end{eqnarray*}
 Using the relation $(A+B)^{\perp}=A^{\perp}\cap B^{\perp}$, which
holds for any subspaces $A$ and $B$ of a symplectic vector space,
one sees that the transversality equation

\begin{eqnarray*}
T_{K(p)}(V\times W)+T_{K(p)}(M\times\Delta_{M}\times P) & = & T_{K(p)}(\overline{M}\times N\times\overline{N}\times P)\end{eqnarray*}
 becomes equivalent to

\begin{eqnarray*}
\underbrace{T_{K(p)}(V\times W)\cap\Big(\{0\}\times T_{(\psi(p),\psi(p))}\Delta_{N}\times\{0\}\Big)}_{U} & = & \{(0,0,0,0)\}.\end{eqnarray*}
 We shall now prove that this last equation holds. By assumption,
we have that

\begin{eqnarray}
T_{(\phi\circ\psi(p),\psi(p))}V\cap\big(T_{\phi\circ\psi(p)}A\times T_{\psi(p)}N\big) & = & T_{(\phi\circ\psi(p),\psi(p))}\graph\phi\label{eq:genmorphV}\\
T_{(\psi(p),p)}W\cap\big(T_{\psi(p)}B\times T_{p}P\big) & = & T_{(\psi(p),p)}\graph\psi.\label{eq:genmorphW}\end{eqnarray}
 Moreover, we may rewrite $U$ as\begin{eqnarray*}
U & = & \underbrace{\bigg(T_{(\phi\circ\psi(p),\psi(p))}V\cap\big(\{0\}\times T_{\psi(p)}N\big)\bigg)}_{J}\times_{T_{\psi(p)}N}\\
 &  & \qquad\underbrace{\times_{T_{\psi(p)}N}\bigg(T_{(\psi(p),p)}W\cap\big(T_{\psi(p)}N\times\{0\}\big)\bigg)}_{L}\end{eqnarray*}
 Equation \eqref{eq:genmorphV} tells us that $J\subset T_{\phi\circ\psi(p)}A\times T_{\psi(p)}B$.
Using this and Equation \eqref{eq:genmorphW}, we can then write

\begin{eqnarray*}
U & = & \bigg(T_{(\phi\circ\psi(p),\psi(p))}\graph\phi\cap\Big(\{0\}\times T_{\psi(p)}B\Big)\bigg)\times_{T_{\psi(p)}N}\\
 &  & \quad\quad\times_{T_{\psi(p)}N}\bigg(T_{(\psi(p),p)}\graph\psi\cap\Big(T_{\psi(p)}B\times\{0\}\Big)\bigg).\end{eqnarray*}
Since \begin{eqnarray*}
T_{(\psi(p),p)}\graph\psi & = & \Big\{(T_{p}\psi(v),v):\; v\in T_{p}C\Big\}\end{eqnarray*}
we see that \begin{eqnarray*}
T_{(\psi(p),p)}\graph\psi\cap\Big(T_{\psi(p)}B\times\{0\}\Big) & = & \{(0,0)\}\end{eqnarray*}
and finally, that $U=\{(0,0,0,0)\}$, as desired. 
\end{proof}
\begin{prop}
The composition of two symplectic micromorphisms $([V],\phi)$ and
$([W],\psi)$ via \eqref{eq:composition} is well defined and yields
a symplectic micromorphism again. 
\end{prop}
\begin{proof}
Lemma \ref{lem:Transvers} together with a continuity argument yield
that there is a neighborhood $U$ of $\graph\phi\times_{B}\graph\psi$
where $V\times W$ and $M\times\Delta_{N}\times P$ still intersect
transversally. Therefore, the map\begin{eqnarray*}
\Red:(V\times_{N}W)\cap U & \longrightarrow & M\times P,\end{eqnarray*}
restricted to this neighborhood, is a immersion according to Theorem
\ref{thm:clean intersection}. At this point, recall that a proper
immersion $i:X\rightarrow Y$ that is injective on a closed submanifold
$A\subset X$ is a embedding on a neighborhood of $A$. Since the
maps $\phi$ and $\psi$ are smooth, $\graph\phi\times_{B}\graph\psi$
is closed. Moreover, on this submanifold, we have that\begin{eqnarray*}
\Red\bigg(\Big(\big(\phi\circ\psi\big)(p),\psi(p),\psi(p),p\Big)\bigg) & = & \Big(\big(\phi\circ\psi\big)(p),p\Big),\end{eqnarray*}
meaning that $\Red$ maps $\graph\phi\times_{B}\graph\psi$ diffeomorphically
to $\graph(\phi\circ\psi)$. Therefore, there is a neighborhood $\overline{U}$
of $\graph\phi\times_{B}\graph\psi$ such that $\Red(\overline{U})$
is a lagrangian submanifold containing $\graph(\phi\circ\psi)$. This
proves that the lagrangian submanifold germ $([W\circ V],\phi\circ\psi)$
is well defined. We need to show that it is a symplectic micromorphism;
i.e. conditions \eqref{eq:CondI} and \eqref{eq:CondII} hold. To
begin with, notice that

\begin{eqnarray*}
W\circ V\cap(A\times P) & = & \Red\Big(\big(V\cap(A\times N)\big)\times_{N}W\Big).\end{eqnarray*}
Now, since, by assumption, $V\cap(A\times N)=\graph\phi\subset A\times B$,
we see that

\begin{eqnarray*}
\big(V\cap(A\times N)\big)\times_{N}W & = & \big(V\cap(A\times N)\big)\times_{N}\big(W\cap(B\times P)\big)\\
 & = & (\graph\phi)\times_{B}(\graph\psi).\end{eqnarray*}
Therefore, we obtain \eqref{eq:CondI} for $W\circ V$, namely \begin{eqnarray*}
W\circ V\cap(A\times P) & = & \Red\Big((\graph\phi)\times_{N}(\graph\psi)\Big)\\
 & = & \graph(\phi\circ\psi).\end{eqnarray*}
Set $K(p)=\big(\phi\circ\psi(p),\psi(p),\psi(p),p\big)$ with $p\in P$.
Realizing that \begin{eqnarray*}
T_{\big(\phi\circ\psi(p),p\big)}(W\circ V) & = & T_{K(p)}\Red\bigg(\Big(T_{\big(\phi\circ\psi(p),\psi(p)\big)}V\Big)\times_{T_{\psi(p)}N}\Big(T_{\big(\psi(p),p\big)}W\Big)\bigg),\end{eqnarray*}
a similar computation on the tangent space level yields \eqref{eq:CondII}
for $W\circ V$. 
\end{proof}

\section{Symplectic categories}

In this section, we reinterpret the results obtained so far in the
language of monoidal categories. We refer the reader to \cite{MacLane1998}
for an exposition on monoidal categories.

\begin{notation*}
We will sometimes write $\mathcal{C}_{0}$ for the objects and $\mathcal{C}_{1}$
for the morphisms of a category $\mathcal{C}$. Accordingly, given
a functor $F:\mathcal{C}\rightarrow\mathcal{D}$, we denote by $F_{0}:\mathcal{C}_{0}\rightarrow\mathcal{D}_{0}$
the object component of $F$ and by $F_{1}:\mathcal{C}_{1}\rightarrow\mathcal{D}_{1}$
its morphism component. 
\end{notation*}
So far, we have seen that {}``symplectic categories'' come in four
flavors. In the macroworld of symplectic manifolds we have: 

\begin{itemize}
\item $\Sympl$, the usual symplectic category of symplectic manifolds and
symplectomorphisms, 
\item $\ExtSympl$, the extended symplectic {}``category'', where symplectomorphisms
are replaced by canonical relations, and which is not a category.
\end{itemize}
In the microworld of symplectic microfolds, we have: 

\begin{itemize}
\item $\MicroSympl$, the $\mathbf{microsymplectic}$ $\mathbf{category}$;
i.e. the category of symplectic microfolds and symplectomorphism germs, 
\item $\ExtMicroSympl$, the $\mathbf{extended}$ $\mathbf{microsymplectic}$
$\mathbf{category}$; i.e., the category of symplectic microfolds
and symplectic micromorphisms. 
\end{itemize}
A major improvement in the microworld is that, this time, symplectic
micromorphisms always compose. Hence, $\ExtMicroSympl$ is a category,
which enlarges the category $\MicroSympl$ in the following precise
sense:

\begin{defn}
\label{def:Extension}A category $\mathcal{D}$ is said to be an $\mathbf{enlargement}$
of a category $\mathcal{C}$ if there is a functor $F:\mathcal{C}\rightarrow\mathcal{D}$
such that $F_{0}$ is a bijection and such that $F_{1}$ is injective
and bijective on the isomorphisms; i.e.\begin{eqnarray*}
\operatorname{Iso}(x,y) & \simeq & \operatorname{Iso}\big(F_{0}(x),F_{0}(y)\big),\end{eqnarray*}
for all objects $x,y\in\mathcal{C}_{0}$. In this case, we also call
the functor an $\mathbf{enlargement}$. 
\end{defn}
Intuitively, enlarging a category means keeping the same objects while
adding morphisms that are not isomorphisms.

Clearly, the functor\begin{eqnarray*}
\GR:\MicroSympl & \longrightarrow & \ExtMicroSympl\end{eqnarray*}
that is the identity on objects and that takes a symplectomorphism
germ to its graph is an enlargement of categories.

The extended microsymplectic category is a symmetric monoidal category.
The tensor product of symplectic microfolds is simply given by\begin{eqnarray*}
[M,A]\otimes[N,B] & := & [M\times N,A\times B].\end{eqnarray*}
Given two symplectic micromorphisms\begin{eqnarray*}
([V_{i}],\phi_{i}):[M_{i},A_{i}] & \longrightarrow & [N_{i},B_{i}],\quad i=1,2,\end{eqnarray*}
we define their tensor product as\begin{eqnarray*}
([V_{1}],\phi_{1})\otimes([V_{2}],\phi_{2}) & := & \bigg(\Big[\big(\id_{M_{1}}\times\epsilon_{N_{1},M_{2}}\times\id_{N_{2}}\big)(V_{1}\times V_{2})\Big],\phi_{1}\times\phi_{2}\bigg),\end{eqnarray*}
where $\epsilon_{X,Y}(x,y)=(y,x)$ is the usual factor permutation.
The unit object $\neutral$ is the cotangent bundle of the one-point
manifold $\{\star\}$. As shown in Section \ref{exa: Unit object and morphism},
$\neutral$ is initial. The symmetry isomorphisms are given by\begin{eqnarray*}
\sigma_{[M,A],[N,B]} & := & \Big(\big[\graph\epsilon_{M,N}\big],\epsilon_{B,A}\Big).\end{eqnarray*}
Note that the opposite symplectic manifold $\overline{(M,\omega)}=(M,-\omega)$
has its natural micro version $\overline{[M,L]}:=[\overline{M},L]$.
It is straightforward, although cumbersome, to verify the following:

\begin{thm}
$(\ExtMicSym,\otimes,E,\sigma)$ is a symmetric monoidal category
with initial unit $E$. 
\end{thm}
We conclude this section by commenting on the relationship between
the extended microsymplectic category and the lagrangian operads introduced
in \cite{Cattaneo,CDF2005}. 

\begin{defn}
An $\mathbf{operad}$ is a collection $\{A(n)\}_{n\geq0}$ of sets
together with composition laws\begin{eqnarray*}
A(n)\times A(k_{1})\times\cdots\times A(k_{n}) & \longrightarrow & A(k_{1}+\cdots+k_{n})\\
(F,G_{1},\dots,G_{n}) & \longmapsto & F(G_{1},\dots,G_{n})\end{eqnarray*}
for each $n,k_{1},\dots,k_{n}\in\N$, satisfying the associativity
equations \begin{align*}
\Big(F(G_{1},\dots,G_{n})\Big)(H_{11},\dots,H_{1k_{1}},\dots H_{n1},\dots,H_{nk_{n}})=\\
F\Big(G_{1}(H_{11},\dots,H_{1k_{1}}),\dots,G_{n}(H_{n1},\dots,H_{nk_{n}})\Big),\end{align*}
and unit $I\in A(1)$ such that $F(I,\dots,I)=F$ for all $F\in A$. 
\end{defn}
For any object $X$ in a monoidal category $(\mathcal{C},\otimes,\neutral)$,
one defines the endomorphism operad $\Endop(X)$ of $X$ to be the
collection\begin{eqnarray*}
\Endop(X)(n) & = & \hom(X^{\otimes n},X)\end{eqnarray*}
with the usual convention that $X^{\otimes n}=\neutral$ for $n=0$.
The composition laws are given by the tensor product and the usual
composition in the category:\begin{eqnarray*}
F(G_{1},\cdots,G_{n}) & := & F\circ(G_{1}\otimes\cdots\otimes G_{n}).\end{eqnarray*}
The unit is the identity morphism $\id_{X}\in\hom(X,X)$. 

Since $\ExtMicroSympl$ is a monoidal category, it makes sense to
consider the endomorphism operad $\Endop([M,A])$ of a symplectic
microfold $[M,A]$. There are two special operads sitting inside of
it. First, the cotangent lifts of the $n$-diagonal maps $\Delta^{n}:A\rightarrow A^{n}$,
$n\geq1$, form an operad\begin{eqnarray*}
\mathcal{L}_{\Delta}([M,A])(n) & := & \big\{\Cot\Delta^{n}\big\},\quad n\geq1,\\
\mathcal{L}_{\Delta}([M,A])(0) & := & \big\{\neutralmorph_{[M,A]}\big\},\end{eqnarray*}
thanks to the properties\begin{eqnarray*}
\Cot\Delta^{1} & = & \id_{[M,A]}\\
\Delta^{k_{1}+\cdots+k_{n}} & = & \Delta^{n}\circ(\Delta^{k_{1}}\times\cdots\times\Delta^{k_{n}})\\
\Cot\Delta^{n-1} & = & \Cot\Delta^{n}\circ(\id_{[M,A]}\otimes\cdots\otimes\neutralmorph_{[M,A]}\otimes\cdots\otimes\id_{[M,A]}).\end{eqnarray*}
Now, $\mathcal{L}_{\Delta}([M,A])$ sits in the suboperad $\mathcal{L}([M,A])$
of $\Endop([M,A])$ defined as follows. For $n\geq1$, $\mathcal{L}([M,A])(n)$
is the set of symplectic micromorphisms $[M,A]^{\otimes n}\rightarrow[M,A]$
whose core map is the $n$-diagonal $\Delta^{n}$. For $n=0$, we
set $\mathcal{L}([M,A])(0)=\{\neutralmorph_{[M,A]}\}$. Note that
the first degree of this suboperad is interesting: $\mathcal{L}([M,A])(1)$
is the group of symplectorphism germs $[\psi]:[M,A]\rightarrow[M,A]$
fixing $A$. 

In \cite{CDF2005}, $\mathcal{L}_{\Delta}([\Cot\R^{n},\R^{n}])$ was
called the cotangent lagrangian operad over $\Cot\R^{n}$ and $\mathcal{L}([\Cot\R^{n},\R^{n}])$
the local lagrangian operad over $\Cot\R^{n}$. They were introduced
ad hoc in terms of generating functions of lagrangian submanifold
germs.

\end{document}